\documentclass[10pt,widea4paper]{article}
\usepackage{t1enc}
\usepackage[latin2]{inputenc} 
\usepackage[dvips,final]{graphicx} 
\usepackage{tabularx} 
\usepackage{amsfonts}
\usepackage{amsthm}
\usepackage{amssymb}
\usepackage{epsf}

\def\du{\unskip\smash{\lower 1.4ex \hbox{\char34}}\kern-.2ex}
\def\hu{\kern-.2ex\hbox{\char92}}

\newcommand{\bdis}{\begin{displaymath}}
\newcommand{\edis}{\end{displaymath}}
\newcommand{\be}{\begin{equation}}
\newcommand{\ee}{\end{equation}}
\newcommand{\ep}{\epsilon}
\def\Diff{\mbox{\it Diff }}
\def\diff{\mbox{\it diff }}
\def\Der{\mbox{\it Der }}
\def\deg{\mbox{\it deg }}
\def\M{\mfrak{M}}
\def\R{\mbb{R}^{0|1}}
\def\X{\mfrak{X}}
\def\C{C^\infty}
\def\ua{\uparrow}
\def\Om{\Omega}
\def\a{\alpha}
\def\b{\beta}
\def\Ber{\mbox{\it{Ber} }}
\def\det{\mbox{\it{det} }}

\newcommand{\ra}{\rightarrow}

\newcommand{\mcal}{\mathcal}
\newcommand{\mfrak}{\mathfrak}
\newcommand{\pd}{\partial}
\newcommand{\mbb}{\mathbb}

\theoremstyle{citing}

\begin{document}
\title{{\textbf{Pseudodifferential forms and supermechanics}}}
\author{Denis Kochan \\
{\it{Department of Theoretical Physics}}, \\ {\it{Faculty of Mathematics Physics and Informatics}},
\\ {\it{ Comenius University}},\\ {\it{ Mlynsk\' a Dolina F2, 842 48 Bratislava, Slovakia}} \\
 \tt{e-mail: kochan@sophia.dtp.fmph.uniba.sk} }
\date{}
\maketitle
{\abstract{\noindent We investigate (pseudo)differential forms in the
framework of supergeometry. Definitions, basic properties and Cartan calculus
(DeRham differential, Lie derivative, inner product, Hodge operator) are
presented; the symplectic supermechanics (even and odd) is formulated;
and the question of quantization is discussed. In the framework of
supermechanics, we investigate also classical Hamiltonian systems
converting to SUSY-QM after quantization.
}}

\section{Introduction}

Supergeometry is surely an interesting and fruitful branch of mathematics
with a variety of powerful applications in modern theoretical physics,
in particular in SUSY, supergravity and superstrings. From a purely mathematical
point of view, supergeometry is a natural extension of the ordinary
differential geometry by Grassmann variables. Such anti-commuting extensions
represent an essential and inspiring feature of all supermathematics.
The first "supermathematician" was undoubtedly Russian mathematician F. A.
Berezin, who formulated the main principles of supermathematics (see his
famous book \cite{Berezin}). In present days there is a lot of books and
articles about supergeometry and its application in physics, but because our
aim is not to describe all historical circumstances, we refer only to a few
of them \cite{Berezin}-\cite{Nelson} (see also references therein).

Differential forms and Cartan calculus are very effective tools of theoretical
physics. It is well known that nine tenths from classical physics
may be formulated and investigated from the geometrical point of view. The
supergeometrical generalization is therefore very useful when describing
systems with both bosonic and fermionic degrees of freedom. The
supermathematics (inspired by the physics) is a right tool to do
this rigorously. The aim of this paper is to present a very short exposition
on pseudodifferential forms (roughly speaking differential forms
on an arbitrary supermanifold) and related Cartan calculus. It will be
shown that such objects are just functions on another (bigger) supermanifold,
therefore they are very easy to handle. All standard geometrical
operations will be encoded into distinguished vector fields whereupon the
calculus will become simple and beautiful.

\section{Pseudodifferential forms}

It is well known (for more details see \cite{Manin}) that the
{\it{\textbf{pseudodifferential forms}}} on
arbitrary smooth $m|n$-dimensional supermanifold $\M$
are defined as functions on $\Pi T\M$ (odd tangent bundle).
The standard differential operations on the forms (DeRham differential,
Lie derivative, inner product) are identified with
special supervector fields on $\Pi T\M$. To obtain their
exact expressions we use the fruitful idea of
Maxim Kontsevich, who pointed out that (see \cite{Kontsevich})
\be
\Pi T\M\equiv\{\mbox{\ supermaps:\ }\R\rightarrow\M\}\ .
\ee
The supergroup $\Diff(\R)$ defines via its natural right action
\bdis
\Pi T\M\times\Diff(\R)\ra\Pi T\M \ \ \ \ \ \ \ \ \ \ \  (F,g)\mapsto F\circ g
\edis
the left invariant (fundamental) vector fields
$Q,\,E\in\X(\Pi T\M)$. Their expression in arbitrary local coordinates
$(x^i,\xi^\a;\psi^i,y^\a)=(x^i,\xi^\a;{\rm{d}}x^i,{\rm{d}}\xi^\a)$
covering the odd tangent bundle $(i=1,\dots ,m \mbox{\ and\ } \a=1,\dots ,n)$
is very simple, namely
\begin{eqnarray}
Q & = & \psi^i\pd_{x^i\,}+y^\a\pd_{\xi^\a} \ \ \ \ \mbox{\it{DeRham differential}}\ ,\label{deRham}\\
E & = & \psi^i\pd_{\psi^i}+y^\a\pd_{y^\a}\ \ \ \ \mbox{\it{Euler field}}\ .\label{Euler}
\end{eqnarray}
Throughout the paper, we use left derivatives with respect to Grassmann
variables and Einstein's summation convention. The Euler vector field "measures"
the degree of homogenity of pseudodifferential forms under the
supergroup action, therefore the superalgebra $\C(\Pi T\M)$ has also a
$\mbb{Z}$-graded structure. A direct calculation gives the (super)commutation
relations in the superalgebra $\diff(\R)$
\begin{eqnarray}
{}[Q,\, Q] & = & 2Q^2=0\ ,\\
{}[E,\, E] & = & 0 \ ,\\
{}[E,\, Q] & = & Q\ .
\end{eqnarray}

Now we shall very briefly sketch the genesis of Lie derivatives and inner
product on pseudodifferential forms under some
vector field $V=V^i(x,\xi)\pd_{x^i}+V^\a(x,\xi)\pd_{\xi^a}=$ $V(x^i)\pd_{x^i}+V(\xi^\a)\pd_{\xi^\a}\in\X(\M)$.
The Lie supergroup $\Diff(\M)$ acts on the supermanifold $\Pi T\M$,
\bdis
\Diff(\M)\times\Pi T\M\ra\Pi T\M \ \ \ \ \ \ \ \ \ \ \  (g,F)\mapsto g\circ F\ ,
\edis
and therefore for arbitrary element $V$ of the superalgebra $\diff(\M)=\X(\M)=\Der(\C(\M))$
we have some vector field $V^\uparrow$ defined on $\Pi T\M$. Notation $V^\uparrow$
reflects the obvious fact that we have lifted the vector field $V$ from the supermanifold $\M$
to $\Pi T\M$. A straightforward coordinate computation\footnote{In the case of an odd vector field
($\tilde{V}=1$) it is necessary to consider the {\emph{superflow}} (homomorphism of the supergroups $\mbb{R}^{1|1}$
and $\Diff(\M)$), whose infinitesimal ($\Delta t,\Delta\ep$) version
in the coordinates is
\bdis
(x^i;\,\xi^\a)\mapsto \bigl(x^i+\Delta\ep V(x^i)+\frac{\Delta t}{2}[V,\, V](x^i);\, \xi^\a+\Delta\ep V(\xi^\a)+\frac{\Delta t}{2}[V,\, V](\xi^\a)\bigr)\ .
\edis} gives
\be\label{Lie derivation}
V^\ua=V(x^i)\pd_{x^i}+V(\xi^\a)\pd_{\xi^\a}+(-1)^{\tilde{V}}Q(V(x^i))\pd_{\psi^i}+(-1)^{\tilde{V}}Q(V(\xi^\a))\pd_{y^\a}\ .
\ee
The procedure of lifting vector fields preserves parity, $\tilde{V^\ua}=\tilde{V}$.

Apart from this natural lifting construction it is also possible to assign to any $V\in\X(\M)$
certain vector field $V_\ua$ on $\Pi T\M$ such that $\tilde{V_\ua}=\tilde{V}+1$ and
\be\label{Cartan magic formula}
{}[V_\ua,\, Q]=V^\ua\ .
\ee
Obviously, the coordinate expression for $V_\ua$ is
\be\label{inner product}
V_\ua=V(x^i)\pd_{\psi^i}+V(\xi^\a)\pd_{y^\a}\ .
\ee
It is easy to confirm the validity of supercommutations relations
\be
\begin{array}{lcccl}
   {} [E,\, V^\ua]\ =0\ , & & & & [E,\, V_\ua]=-V_\ua\ , \\
   {} [V^\ua,\, Q]\ =0\ , & & & & [V^\ua,\, W^\ua]=[V,\, W]^\ua \ ,\\
  {}[V_\ua,\, W_\ua]=0\ , & & & & [V^\ua,\, W_\ua]\,=[V,\, W]_\ua\ .
\end{array}
\ee
The vector fields $V^\ua$
corresponds to the {\emph{Lie derivative}} $\mcal{L}_V$ (with respect to $V$),
whereas $V_\ua$ corresponds to the {\emph{inner product}} $i_V$ (with $V$). The
equation (\ref{Cartan magic formula}) is the famous {\emph{Cartan formula}}.

An arbitrary pseudodifferential form is a function on the supermanifold $\Pi T\M$ and therefore
it may be expressed in any local coordinates as
\be
f=f(x,\xi,\psi,y)=\sum\limits_{{\beta\geq 0 \atop k\geq 0}}\sum\limits_{{\a_1,\dots ,\a_\beta \atop i_1,\dots ,i_k}}f_{\a_1,\dots ,\a_\beta,i_1,\dots ,i_k}(x,y)\xi^{\a_1}\dots\xi^{\a_\beta}\psi^{i_1}\dots\psi^{i_k}\ ,
\ee
where the ordinary real (complex) valued functions $f_{\a_1,\dots ,i_k}(x,y)$
are skew-sym\-metric in
the indices ${\a_1,\dots ,i_k}$. The Berezin integral (for more details see
\cite{Berezin},\cite{Manin}) of a function $f$ on $\Pi T\M$,
\be\label{integral}
I[f]:=\int\limits_{\M}\overline{{\rm d}x{\rm d}\xi}\int\limits_{\mbb{R}^{n|m}}\overline{{\rm d}\psi{\rm d}y}\ f(x,\xi,\psi, y),
\ee
defines the {\emph{integral}} of the pseudodifferential form $f$ over $\M$. It is clear that such integral
is not well defined for all elements of $\C(\Pi T\M)$, because the supermanifold $\mbb{R}^{n|m}$ (the typical
fiber in the bundle $\Pi T\M\ra\M$) is not compact. The coordinate transformation on the supermanifold $\M$
$(x^i,\xi^\a)\mapsto \bigl(X^i(x,\xi),\Xi^\a(x,\xi)\bigr)$ induces the corresponding transformation on
$\Pi T\M$,
\bdis
(x^i,\xi^a;\psi^i,y^\a)\mapsto \biggl(X^i,\Xi^\a;\Psi^i=\psi^k\frac{\pd X^i}{\pd x^k}+y^\beta\frac{\pd X^i}{\pd\xi^\beta},Y^\a=\psi^k\frac{\pd\Xi^\a}{\pd x^k}+y^\beta\frac{\pd \Xi^\a}{\pd\xi^\beta}\biggr)\ .
\edis
Its berezinian is equal to unity, therefore the integral defined in (\ref{integral}) is coordinate independent.

It is also possible to define (at least formally) Hodge $*$ operator acting on
pseudodifferential forms. As in the case of ordinary differential forms,
the essential ingredient to define Hodge $*$ operator is metric. The metric may
be introduced on an arbitrary $m|n$-dimensional smooth supermanifold $\M$
(in particular, on an ordinary manifold $M$) as an even regular
(non-degenerate) quadratic function at the tangent bundle $T\M$
\be\label{metric}
G(x,\xi,z,\sigma)=(z^i,\sigma^\a)
\left(
\begin{array}{cc}
g_{ij}(x,\xi) & \Gamma_{i\beta}(x,\xi) \\
\Gamma_{\a j}(x,\xi) & h_{\a\beta}(x,\xi)
\end{array}
\right)
\left(
\begin{array}{c}
z^j \\
\sigma^\beta
\end{array}
\right)\ .
\ee
The non-degeneracy condition reads
\be
\Ber G=\det(g_{ij}-\Gamma_{i\beta}(h_{\a\beta})^{-1}\Gamma_{\a j})\det(h_{\a\beta})^{-1}\neq 0\ .
\ee
The transformation of local coordinates on the supermanifold $\M$ and
the transformation of the coordinates $(x^i,\xi^\a,z^i,\sigma^\a)$ on $T\M$
are coupled by
\bdis
(x^i,\xi^a;z^i,\sigma^\a)\mapsto \biggl(X^i,\Xi^\a;Z^i=z^k\frac{\pd X^i}{\pd x^k}+\sigma^\beta\frac{\pd X^i}{\pd\xi^\beta},\Sigma^\a=z^k\frac{\pd\Xi^\a}{\pd x^k}+\sigma^\beta\frac{\pd \Xi^\a}{\pd\xi^\beta}\biggr)\ .
\edis
Let us emphasize that non-degeneracy of $G$ implies that the even skew-symmetric
matrix $h_{\a\beta}$ is invertible, consequently, $m|n$-dimensional supermanifold $\M$ may be Riemannian
only if $n$ is even (this fact is strongly reminiscent of the situation in symplectic geometry).
For the functions on $\Pi T\M$ (pseudodifferential forms on $\M$) the {\emph{Hodge star operator}} is defined via
Fourier transformation in the "fibre variables" $\psi^i$ and $y^\a$, namely
\be\label{Hodge *}
(*f)(x,\xi,\psi,y):=\int\limits_{\mbb{R}^{n|m}}\frac{\overline{{\rm d}\psi^\prime {\rm d}y^\prime}}{\sqrt{\Ber G}} f(x,\xi,\psi^\prime,y^\prime)
\,e^{-\imath\langle\psi^\prime,y^\prime|\psi,y\rangle_G}\ ,
\ee
where the symbol $\langle\psi^\prime,y^\prime|\psi,y\rangle_G$ denotes the "scalar product" with respect
to the metric $G$
\be\label{scalar product}
\langle\psi^\prime,y^\prime|\psi,y\rangle_G:=({\psi^\prime}^i,{y^\prime}^\a)
\left(\begin{array}{cc} g_{ij} & \Gamma_{i\beta} \\
                 \Gamma_{\a j} & h_{\a\beta}
\end{array}\right)\left(
\begin{array}{c}
{-\psi^j} \\
{y^\beta}
\end{array}
\right)\ .
\ee
A straightforward calculation shows that the definition of Hodge $*$ operator
is independent on the choice of coordinates, and for ordinary forms it gives
a multiple of standard $*_g$. As above, the non-compactness of the fibre
$\mbb{R}^{n|m}$ implies that $*$ is defined only for functions that are
behaving well in the variables $y^\a$ at infinity, for example for the
functions with compact support.

\section{Fr\" olicher-Nijenhuis brackets}

In this section we will consider only an ordinary (real) smooth
$m$-dimensional manifold $M$. It has been shown that exterior
algebra of differential forms $\Omega(M)$ is in one to one
correspondence with $\C(\Pi TM)$. The supergroup $\Diff(\R)$
defines via its action on $\C(\Pi TM)$ also a $\mbb{Z}$-graded
structure on the Lie superalgebra $\X(\Pi TM)=\Der(\C(\Pi TM))$,
namely \be V\in\X^{(k)}(\Pi TM)\Leftrightarrow
[E,\,V]=k\,V=:\deg(V)\,V\ . \ee The vector fields from
$\X^{(k)}(\Pi TM)$ encode, from the ordinary point of view, the
derivations of $\Omega(M)$ of degree $k$, and moreover, from the
(super)Jacobi identity it is clear that
$\deg([V,\,W])=\deg(V)+\deg(W)$. The derivations of $\Omega(M)$
commuting with DeRham deferential are represented by special
vector fields from $\X(\Pi TM)$ that commute with $Q$. The Lie
subsuperalgebra of such derivations will be denoted $\X_Q(\Pi TM)$
(the (super)Jacobi identity states that supercommutator of two
$Q$-invariant vector fields is again $Q$-invariant). The obvious
coordinate expression for $V_A\in\X_Q^{(k)}(\Pi TM)$ is \be
V_A=A^i(x,\psi)\pd_{x^i}+(-1)^kQ(A^i(x,\psi))\pd_{\psi^i}\ , \ee
where $A^i(x,\psi)$ are functions on $\Pi TM$ with the degree of
homogenity $k$ (differential forms of $k$th degree). Using
(\ref{Cartan magic formula}) it is possible to assign to any
$Q$-invariant vector field $V_A$ of degree $k$ another vector
field $v_A\in\X^{(k-1)}(\Pi TM)$ (roughly speaking "potential of
$V_A$"), such that \be v_A=A^i(x,\psi)\pd_{\psi^i}=A^i_{j_1,\dots
,j_k}(x)\psi^{j_1}\dots\psi^{j_k}\pd_{\psi^i}\ . \ee There is a
unique correspondence between $\X_Q^{(k)}(\Pi TM)$ and total
skew-symmetric tensor fields of type ${1 \choose k}$ on manifold
$M$, because any such tensor $A$ is completely characterized by a
set of component functions $A^i_{j_1,\dots ,j_k}(x)$. We define
the \emph{Fr\" olicher-Nijenhuis brackets} between two such
tensors $A,\,B$ of ranks ${1 \choose k},\,{1 \choose l}$
respectively by the
 supercommutator of corresponding vector fields, namely
\be\label{Fr\" olicher-Nijenhuis brackets}
[V_A,\,V_B]=:V_{\{A,\,B\}_{NB}}\ . \ee All properties of the Fr\"
olicher-Nijenhuis brackets may be obtained from the defining
equation (\ref{Fr\" olicher-Nijenhuis brackets}) and from the
properties of the supercommutator (for more details see Nijenhuis
pioneering work \cite{Nijenhuis}).

\section{Symplectic supermechanics}

The Cartan calculus is undoubtedly a useful tool in modern theoretical physics. A very nice
and simple example of its application is the symplectic formulation of Hamiltonian mechanics
(see for example \cite{Arnold},\cite{Abraham-Marsden}). Our next task is
formulation and brief description of the supersymmetric extension (via pseudodifferential
forms) of ordinary Hamiltonian mechanics, and quantization of the extended theory.
In the general Poisson setting this has been done in \cite{Martin}-\cite{Casalbuoni2}.

The supermanifold $\M$, endowed with a particular function (pseudodifferential form)
$\Om\in\C(\Pi T\M)$, is called\vspace{0,3cm}

\begin{minipage}{7,6cm}
{\it{even-symplectic}}:\\
if $\mbox{\it{dim}}(\M)=2m|n$\ ;\\
$\tilde{\Om}=0$\ ;\\
$Q(\Om)=0$ (closedness)\ ;\\
$E(\Om)=2\Om$ ($\Om$ is 2-form) $\Rightarrow $\\
$\Om$ is regular polynomial \\
(non-degenerate) of $2$-nd \\
degree in $\psi, y$\ .
\end{minipage}\ \
\begin{minipage}{7,6cm}
{\it{odd-symplectic}}:\\
if $\mbox{\it{dim}}(\M)=m|m$\ ;\\
$\tilde{\Om}=1$\ ;\\
$Q(\Om)=0$ (closedness)\ ;\\
$E(\Om)=2\Om$ ($\Om$ is 2-form) $\Rightarrow $\\
$\Om$ is regular polynomial \\
(non-degenerate) of $2$-nd \\
degree in $\psi, y$\ .
\end{minipage}\vspace{0,3cm}
The presence of a symplectic structure on the supermanifold $\M$ allows us to define
Poisson brackets on $\C(\M)$: for arbitrary homogeneous function (in the sense of parity)
$f\in\C(\M)\subset\C(\Pi T\M)$, we define the Hamiltonian vector field
$\zeta_f\in\X(\M)$ by the condition
\be
{\zeta_f}_\ua\,\Om=-(-1)^{\tilde{f}}Q(f)\ .
\ee
It is evident that $\tilde{\zeta_f}=\tilde{f}+\tilde{\Om}$ and the assignment $f\mapsto\zeta_f$
is $\mbb{R}$-linear. The corresponding {\it{Poisson brackets}} are
\be
\{f,\,g\}_{\tilde{\Om}}:={\zeta_f}_\ua {\zeta_g}_\ua\,\Om(-1)^{\tilde{\Om}+\tilde{g}}=\zeta_f^\ua\, g(-1)^{\tilde{\Om}+1}=\zeta_f \,g(-1)^{\tilde{\Om}+1}\ .
\ee
Following formulas are valid
\begin{eqnarray}
\widetilde{\{f,\, g\}_{\tilde{\Om}}} & = & \tilde{f}+\tilde{g}+\tilde{\Om}\ ,\\
\zeta_{\{f,\, g\}_{\tilde{\Om}}} & = & [\zeta_f,\, \zeta_g]\ ,\\
\{f,\, g+h\}_{\tilde{\Om}} & = & \{f,\, g\}_{\tilde{\Om}} + \{f,\, h\}_{\tilde{\Om}}\label{linearity}\ ,\\
\{f,\, g\}_{\tilde{\Om}} & = & -\{g,\, f\}_{\tilde{\Om}}(-1)^{(\tilde{f}+\tilde{\Om})(\tilde{g}+\tilde{\Om})}\label{graded skew-symmetry}\ ,\\
\{f,\, gh\}_{\tilde{\Om}} & = & \{f,\, g\}_{\tilde{\Om}}h + g\{f,\, h\}_{\tilde{\Om}}(-1)^{(\tilde{f}+\tilde{\Om})\tilde{g}}\label{graded Leibniz rule}\ ,\\
0 & = & \{f,\, \{g,\, h\}_{\tilde{\Om}}\}_{\tilde{\Om}}(-1)^{(\tilde{f}+\tilde{\Om})(\tilde{h}+\tilde{\Om})} + \mbox{cyclic permutations}\label{graded Jacobi identity}\ .
\end{eqnarray}
Observables of supermechanics are by definition functions on $\M$; $\C(\M)$ is $\mbb{Z}_2$-graded
commutative algebra and at the same time Lie superalgebra, and the compatibility of these two structures is
guaranteed by the graded Leibniz rule (\ref{graded Leibniz rule}).

To define the dynamics we need the
Hamiltonian $H\in\C(\M)$ (homogeneous element in the sense of parity). The (super)time evolution
is generated by (super)flow $\Phi_{(t,\ep)}(\zeta_H)$ of the corresponding Hamiltonian vector field $\zeta_H$;
the general formula for the pull-back of the (super)flow on the observables is
\be\label{Hamiltonian evolution superflow}
f\mapsto f_{(t,\ep)}:=\biggl\{\begin{array}{ll}e^{-t\zeta_H^\ua}f & \mbox{\ \ when\ \ } \tilde{H}+\tilde{\Om}=0 \ ,\\
e^{-\ep\zeta_H^\ua-\frac{t}{2}\zeta_{\{H,H\}_{\tilde{\Om}}}^\ua}f & \mbox{\ \ when\ \ }\tilde{H}+\tilde{\Om}=1\ .
\end{array}\biggr.
\ee
It is convenient to rewrite the integral formula for the
(super)time\footnote{We see that superflow forces us
to extend the ordinary evolution of classical mechanics in $t\in\mbb{R}$ (additive Lie group)
to superevolution in $(t,\ep)\in\mbb{R}^{1|1}$ (one of the simplest Lie supergoups).
The multiplication in $\mbb{R}^{1|1}$ is given by
$(t,\ep)\cdot(t^\prime,\ep^\prime)=(t+t^\prime+\ep\ep^\prime,\ep+\ep^\prime)$, with
the neutral element ${\rm e}=(0,0)$ and the inverse element $(t,\ep)^{-1}=(-t,-\ep)$.
Lie superalgebra $\mathbf{r}^{1|1}$ is generated by vector fields $V_0=\pd_t$ and
$V_1=\pd_\ep+\ep\pd_t$ obeying $[V_i,\,V_j]=2ijV_0$.}
evolution (\ref{Hamiltonian evolution superflow})
into a differential one. This amounts to writing down the Hamiltonian
equations of motion
\begin{eqnarray}
\mbox{if\ } \tilde{H}+\tilde{\Om}=0 \mbox{\ \ then\ \ } \{H,\, f\}_{\tilde{\Om}}(-1)^{\tilde{\Om}} & = & \pd_t f\ ,\label{Hamilton1}\\
\mbox{if\ } \tilde{H}+\tilde{\Om}=1 \mbox{\ \ then\ \ } \{H,\, f\}_{\tilde{\Om}}(-1)^{\tilde{\Om}} & = & (\pd_\ep+\ep\pd_t) f\label{Hamilton2}\ .
\end{eqnarray}
The question of symmetry of the Hamiltonian system $(\M,\,\Om,\,H)$ is also very simple:
vector field $V\in\X(\M)$ is a Cartan symmetry if $V^\ua\,\Om=0=V^\ua\, H$. If moreover there
exists some function $F\in\C(\M)$ for which $V_\ua\,\Om=-(-1)^{\tilde{F}}Q(F)$, then we call $V$
an exact Cartan symmetry and the function $F$ is a conserved quantity.

The quantization with a given symplectic structure may be performed by the famous Fedosov construction
\cite{Fedosov}, well known from deformation quantization theory, which for the simplest case
(flat phase space, discussed in more detail below) coincides with the Wigner-Moyal-Weyl
quantization (for more detail see \cite{Flato}).\\
\textbf{Quantization of even supersymplectic structures:} We shall sketch shortly
the quantization procedure for supermanifold $\M=\mbb{R}^{2m|n}$ (equipped with global coordinates
$(x^1,\dots ,x^m,x^{m+1}=p_1,\dots ,x^{2m}=p_m;\xi^1,\dots ,\xi^n)$), the superanalog of ordinary
phase space, endowed with the symplectic structure (pseudodifferential form)
\be\label{canonical form of even symplectic form}
\Om=\psi^i\pi_i-\frac{1}{2}g_{\a\b}y^\a y^\b={\rm d}x^i\wedge{\rm d}p_i-\frac{1}{2}g_{\a\b}{\rm d}\xi^\a\wedge{\rm d}\xi^\b\ ,
\ee
where $g_{\a\b}=\mbox{diag}(+1,\dots ,+1,-1,\dots ,-1)$. This is a natural
extension of the canonical $2$-form $\omega={\rm d}x^i\wedge{\rm d}p_i$ on
$\mbb{R}^{2m}$. It is not so difficult to prove that
the famous Darboux theorem, well known from standard symplectic geometry (see for example \cite{da Silva}),
is valid in supersymplectic case, too, so that an arbitrary $2m|n$-dimensional even symplectic supermanifold $\M$
is locally isomorphic to $\mbb{R}^{2m|n}$ with the symplectic form (\ref{canonical form of even
symplectic form}). Let us stress
that even "supersymplecticity" (contrary to "supermetricity") does not lead to any
restriction on $n$. The corresponding coordinate Hamiltonian vector fields
\be
\zeta_{x^i}=\pd_{p_i}\,,\ \ \ \ \zeta_{p_i}=-\pd_{x^i}\,,\ \ \ \ \zeta_{\xi^\a}=-g^{\a\b}\pd_{\xi^\b}\,,
\ee
determine elementary Poisson brackets
\be\label{commutation relations I}
\begin{array}{lll}
  \{x^i,\, x^j\}=0 &  \{p_i,\, p_j\}=0  & \{p_i,\, x^j\}=\delta_i^j\ , \\
\{x^i,\, \xi^\a\}=0 & \{p_i,\, \xi^\a\}=0 & \{\xi^a,\, \xi^\b\}=g^{\a\b}\ .
\end{array}
\ee
The canonical operator quantization is equivalent, via Weyl isomorphism, to
quantization with operator symbols (ordinary functions on the phase space) and vice versa,
with the operator product replaced by the star product\footnote{which
defines the so called quantum deformation of the classical space of observables}
\begin{eqnarray}
f\star h & = & f h+\frac{\imath\hbar}{2}\{f,\,h\}+o(\hbar) \nonumber \\
& = & f h+\frac{\imath\hbar}{2}\biggl(\frac{\pd f}{\pd
p_i}\frac{\pd h}{\pd x^i}-\frac{\pd f}{\pd x^i}\frac{\pd h}{\pd
p_i}+(-1)^{\tilde{f}+1}\frac{\pd f}{\pd \xi^\a}g^{\a\b}\frac{\pd
h}{\pd \xi^\b}\biggr)+o(\hbar)\ ,
\end{eqnarray}
and the supercommutator $[\hat{f},\,\hat{h}]_\mp$ is in one to one correspondence with a star (super)bracket
$\{f,\, h\}_{\star_\mp}:=f\star h \mp h\star f$. The "star product technique" leads to the canonical supercommutation
relations
\be\label{commutation relations II}
\begin{array}{lll}
{}   [\hat{x}^i,\, \hat{x}^j]_{-\,}=0\ , &    [\hat{p}_i,\, \hat{p}_j]_{-\,}=0\ ,  & [\hat{p}_i,\, \hat{x}^j\,]_{-\,}=\imath\hbar\delta_i^j \ ,\\
{}[\hat{x}^i,\, \hat{\xi}^\a]_{-}=0\ , & [\hat{p}_i,\, \hat{\xi}^\a]_{-}=0\ ,  & [\hat{\xi}^\a,\, \hat{\xi}^\beta]_{+}=\imath\hbar g^{\a \b}\ .
\end{array}
\ee
The appropriate Hilbert space for such QM is $L^2(\mbb{R}^m$, ${\rm d}\mu)\otimes\mbb{C}^{d}$, $d=2^{[\frac{n}{2}]}$,
and the dynamical equations (in Heisenberg picture) coincide with (\ref{Hamilton1}),
(\ref{Hamilton2}) if we put $\tilde{\Om}=0$ and replace $\{.\,,\,.\}_{\tilde{\Om}}$ by $[.\,,\,.]_\mp$.
Now we are ready to present simple examples.\\
{\emph{SUSY-QM}}: In order to simplify the exposition we put
$g^{\a\b}=\delta^{\a\b}$, $m=1$ and $n=2$ (it is well known from the theory of
Clifford algebras \cite{Snygg} that for all even $n$ and all metrics $g^{\a\b}$
the properties of the algebras are similar). The supercommutation
relations (\ref{commutation relations II}) are represented in the Hilbert space
$L^2({\mbb{R}},\,{\rm d}x)\otimes\mbb{C}^2$ by the operators
\be\label{realization of CM II}
\hat{x}=x\,,\ \
\hat{p}=\imath\hbar\frac{{\rm d}}{{\rm d}x}\,,\ \
\hat{\xi}^1=\sqrt{\frac{\hbar}{2}}\left(
\begin{array}{cc}
 0 & \imath \\
 1 & 0
\end{array}\right)\,,\ \
\hat{\xi}^2=\sqrt{\frac{\hbar}{2}}\left(
\begin{array}{cc}
 0 & 1 \\
 \imath & 0
\end{array}\right)\,.
\ee
It is convenient to form from self-adjoint operators $\hat{x}$ and $\hat{p}$ the bosonic
ladder operators
\be
b=\frac{1}{\sqrt{2\hbar}}(\hat{x}-\imath\hat{p})\,,\ \ \ \ \
b^\dag=\frac{1}{\sqrt{2\hbar}}(\hat{x}+\imath\hat{p})\ ,
\ee
and analogically, from anti-self-adjoint operators $\hat{\xi}^\a=\imath(\hat{\xi}^\a)^\dag$
($\a=1,2$) the fermionic ladder operators
\be
f=\frac{1}{\sqrt{2\hbar}}(\hat{\xi}^1-\imath\hat{\xi}^2)\,,\ \ \ \ \
f^\dag=\frac{1}{\sqrt{2\hbar}}(\hat{\xi}^2-\imath\hat{\xi}^1)\ .
\ee
The general even classical Hamiltonian quadratic in the momentum $p$
is of the form
\be
H=\frac{1}{2}(p^2+V_1^2(x))+V_2(x)\xi^2\xi^1\ .
\ee
The most interesting case is $V_2(x)=\pm V_1^\prime(x)$,
because then there exist odd conserved quantities ($\{H,\,q\}=0$)
\be
q_{1\pm}=p\xi^1\mp V_1(x)\xi^2\,,\ \ \ \ \ q_{2\pm}=p\xi^2\mp V_1(x)\xi^1\,,
\ee
such that
\be
\{q_{i+},\,q_{j+}\}=2\delta_{ij}H_+\,,\ \ \ \ \ \{q_{i-},\,q_{j-}\}=2\delta_{ij}H_-\,,
\ee
where
\bdis
H_\pm=\frac{1}{2}(p^2+V_1^2(x))\pm V_1^\prime(x)\xi^2\xi^1\ .
\edis
The corresponding Hamiltonian vector fields $\zeta_q$ generate a special kind
of internal (exact Cartan) symmetry, which turns into supersymmetry
on quantum level. Quantization procedure leads to the Hamiltonian
\be
\hat{H}_{\pm}=\frac{1}{2}\biggl(\hat{p}^2+V_1^2(\hat{x})\pm\hbar V_1^\prime(\hat{x})\sigma_3\biggr)\ ,
\ee
which is the famous Witten's Hamiltonian for the SUSY-QM (proposed in \cite{Witten1},\cite{Witten2} and
studied in \cite{Salomonson},\cite{Crombrugghe}). The operators
$\hat{Q}_{i\pm}=\frac{1}{\sqrt{\imath\hbar}}\hat{q}_{i\pm}$ are odd generators of quantum supersymmetry
\be
[\hat{Q}_{i+},\,\hat{Q}_{j+}]_+=2\delta_{ij}\hat{H}_+\,,\ \ \ \ \ [\hat{Q}_{i-},\,\hat{Q}_{j-}]_+=2\delta_{ij}\hat{H}_-\,.
\ee
Very nice and pedagogical articles about SUSY-QM may be found in \cite{Gendenshtein},
\cite{Grosse} (see also \cite{Ravndal},
where SUSY-QM is built from classical mechanics in a Lagrangian framework).\\
{\emph{{Spin $\frac{1}{2}$ particle}}:
Another important example is $g^{\a\b}=\delta^{\a\b}$, $m=3$ and $n=3$. The operator
realization of (\ref{commutation relations II})
is for even variables the same as in (\ref{realization of CM II}), while odd
variables are represented by the operators
\be
\hat{\xi}^\a=\sqrt{\frac{\imath\hbar}{2}}\sigma^\a \ \ \ \a=1,2,3\ ,
\ee
where $\sigma^\a$ are Pauli matrices (for a general odd $n$ the same is true
with the $\sigma$'s replaced by the generators of corresponding Clifford algebra).
Classical even Hamiltonian
\be\label{Pauli Hamiltonian}
H=\frac{1}{2m}(\vec{p}+q\vec{A})^{2}+q\varphi+\tilde{g}\frac{q}{2m}\frac{\pd A_\b}{\pd x^\a}\xi^\a\xi^\b\ ,
\ee
corresponds after quantization to the well known Pauli Hamiltonian, which describes
non-relativistic spin $\frac{1}{2}$ particle with charge $q$, mass $m$ and Land\' e $\tilde{g}$-factor
in an external electromagnetic field
(more details about quantum mechanics related to Pauli equation could be found e. g. in \cite{Landau}).
For the Hamiltonian (\ref{Pauli Hamiltonian}) and a purely magneto-static field
described by the vector potential $\vec{A}$ it is possible to find a
classical odd conserved quantity
\be
Q=\frac{1}{\sqrt{m}}(p_\mu+qA_\mu)\biggl(\xi^\mu+q(\tilde{g}-2)\,\frac{\ep^{\mu\a\b}\pd_{x^\a}A_\b}{2(\vec{p}+q\vec{A})^{2}}\,\xi^1\xi^2\xi^3\biggr)\ ,
\ee
which obeys $\{Q,\,Q\}=2H$. The corresponding Hamiltonian vector field
$\zeta_Q$ is an exact Cartan symmetry (more details about supersymmetry related to
Pauli equation in specific configurations of magnetic field see in \cite{Niederle}).\\
\textbf{Quantization(?) of odd supersymplectic structures:} The Darboux theorem for odd symplectic
structures (for more details and a proof see \cite{Shander}, \cite{Khudaverdian I}) states that
an arbitrary $m|m$-dimen\-sional odd symplectic supermanifold $\M$ is locally
isomorphic to $\mbb{R}^{m|m}$ (described by global coordinates $(x^1,\dots ,x^m,\xi^1,\dots ,\xi^m)$)
with the odd symplectic structure
\be
\Om=\delta_{i\a}\psi^i y^\a=\delta_{i\a}{\rm d}x^i\wedge{\rm d}\xi^\a\ .
\ee
The coordinate Hamiltonian vector fields
\be
\zeta_{x^i}=-\delta^{i\a}\pd_{\xi^\a}\,,\ \ \ \ \ \ \ \zeta_{\xi^\a}=\delta^{\a i}\pd_{x^i}\,,
\ee
define elementary Poisson brackets (one usually refers to the Poisson brackets for
odd symplectic structures as Buttin brackets)
\be
\{x^i,\, x^j\}=0\,,\ \ \ \ \{\xi^a,\, \xi^\b\}=0\,,\ \ \ \ \{\xi^\a,\, x^i\}=\delta^{\a i}\,.
\ee
The well known and simplest example of an odd symplectic (Poisson) supermanifold is $\Pi T^*M$
(the odd cotangent bundle associated to $T^*M$ by changing parity in fibres).
Functions on $\Pi T^*M$ correspond to multivector fields on the manifold $M$ and Poisson
bracket of two functions is defined as Schouten bracket of congruent multivector fields.

It is not fully clear (at least for us) what exactly means the "quantization" of odd
symplectic structures. It is possible to adopt the quantization schema from the even case,
deforming the algebra $\C(\mbb{R}^{m|m})$ by introducing the associative $\star$ product
\be
f\star g:=fg+\{f,\,\kappa g\}=fg+\kappa\biggl(\frac{\pd f}{\pd \xi^\a}\delta^{\a i}\frac{\pd g}{\pd x^i}+(-1)^{\tilde{f}}\frac{\pd f}{\pd x^i}\delta^{i\a}\frac{\pd g}{\pd \xi^\a}\biggr)\ ,
\ee
with odd deforming parameter $\kappa$ in order to preserve parity.
Unfortunately, such associative star product is $\mbb{Z}_2$-graded commutative and it may be
proved (see \cite{Tamarkin}) that the $\star$ product and the ordinary (exterior) product are
equivalent\footnote{another possible location of $\kappa$ leads to a different supercommutative product
or to a $\star$ product that is not associative}.
It was shown by P. \v Severa (for more details see \cite{Severa}) that an odd Poisson structure
(in particular, an odd symplectic structure) on arbitrary $m|m$-dimensional smooth supermanifold $\M$
may be used to deform/quantize the algebra of pseudodifferential forms
over $\M$.

Let us stress, finally, that odd symplectic structure is a crucial ingredient of Batalin-Vilko\-vis\-ky
formalism (for more details see \cite{Batalin-Vilkovisky I},\cite{Batalin-Vilkovisky II}),
whose geometrical background is discussed in
\cite{Khudaverdian I},\cite{Khudaverdian II}-\cite{Schwarz}.\\
\begin{center} {\textbf{Acknowledgement}} \end{center}

\noindent I would like to express my deep gratitude to my teacher P. \v Severa,
who was an excellent guide of mine to supergeometry and who is the
originator of many ideas in this article. I am also very grateful to P.
Pre\v snajder and V. Balek who gave me much advice, many thanks go to my girlfriend
Tulka ($\heartsuit$) for her mental encouragement.

\end{document}